\documentclass[12pt,a4paper]{article}
\usepackage[utf8]{inputenc}
\usepackage[T1]{fontenc}
\usepackage{amsmath}
\usepackage{amsfonts}
\usepackage{amssymb,amsthm}
\usepackage{graphicx}
\usepackage{mathrsfs}
\usepackage{hyperref}
\usepackage{color}
\usepackage{etoolbox}
\usepackage{tikz}
\usepackage{subcaption} 

\newtheorem{thm}{Theorem}[section]

\newtheorem{pro}[thm]{Proposition}
\newtheorem{cor}[thm]{Corollary}
\newtheorem{lem}[thm]{Lemma}
\newtheorem{rem}[thm]{Remark}

\newcommand{\R}{\mathbb{R}}
\newcommand{\N}{\mathbb{N}}

\def\qq#1{\qquad \mbox{#1}\quad}
\newcommand{\al}{\alpha}
\newcommand{\be}{\beta}

\newcommand{\e}{\varepsilon}

\newcommand{\la}{\lambda}
\newcommand{\La}{\Lambda}

\newcommand{\na}{\nabla}

\newcommand{\Om}{\Omega}
\newcommand{\Omb}{\overline{\Om}}
\newcommand{\p}{\partial}
\newcommand{\s}{\sigma}

\newcommand{\z}{\zeta}
\makeatletter
\patchcmd{\maketitle}{\@fnsymbol}{\@alph}{}{}  
\makeatother
\title{
	Bifurcation and Multiplicity Results for  Elliptic Problems with  Subcritical Nonlinearity on the Boundary }

\makeatletter
\newcounter{author}
\renewcommand*\author[1]{%
	\stepcounter{author}%
	\ifnum\c@author=1
	\gdef\@author{#1}%
	\else
	\xdef\@author{\unexpanded\expandafter{\@author\and#1}}%
	\fi
	\csgdef{author@\the\c@author}{#1}}
\newcommand*\email[1]{%
	\csgdef{email@\the\c@author}{#1}}
\newcommand*\address[1]{%
	\csgdef{address@\the\c@author}{#1}}
\AtEndDocument{%
	\xdef\author@count{\the\c@author}%
	\c@author=1
	\print@authors}
\newcommand*\print@authors{%
	\ifnum\c@author>\author@count
	\else
	\print@author{\the\c@author}%
	\advance\c@author by 1
	\expandafter\print@authors
	\fi}
\newcommand*\print@author[1]{%
	\par\medskip
	\begin{tabular}{@{}l@{}}%
		\textsc{Addresses of \csuse{author@#1}}\\
		\csuse{address@#1}\\
		\textit{E-mail address}:
		\href{mailto:\csuse{email@#1}}{\csuse{email@#1}}
\end{tabular}}
\makeatother
{\small
	\author{Shalmali Bandyopadhyay}
	\address{UNC Greensboro, Greensboro, NC, USA}
	\email{s\_bandyo@uncg.edu}
	
	\author{Maya Chhetri}
	\address{UNC Greensboro, Greensboro, NC, USA}
	\email{m\_chhetr@uncg.edu}
	
	\author{Briceyda B. Delgado}
	\address{Universidad Autónoma de Aguascalientes, Aguascalientes, Mexico}
	\email{bdelgado@math.cinvestav.mx}
	
	\author{Nsoki Mavinga}
	\address{Swarthmore College, Swarthmore, PA, USA}
	\email{nmaving1@swarthmore.edu}
	
	\author{Rosa Pardo}
	\address{Universidad Complutense de Madrid, Madrid, Spain}
	\email{rpardo@ucm.es}
}

\begin{document}
	\maketitle
	{ }\makeatletter{\renewcommand*{\@makefnmark}{}
		\footnotetext{\emph{Keywords:} elliptic problem, nonlinear boundary conditions, superlinear and subcritical, local bifurcation, degree theory,  global bifurcation. }} 
	
	{ }\makeatletter{\renewcommand*{\@makefnmark}{}
		\footnotetext{\emph{Mathematics Subject Classification (2020):} 35J65, 35J61, 35J15.}}

	\begin{abstract}
		We consider an elliptic problem with nonlinear boundary condition involving nonlinearity with superlinear and subcritical growth at infinity and a bifurcation parameter as a factor. We use re-scaling method, degree theory and continuation theorem to prove that 
		there exists a connected branch of positive solutions bifurcating from infinity when the parameter goes to zero. 
		Moreover,  if the nonlinearity satisfies additional conditions near zero, we establish a global bifurcation result, and discuss the number of positive solution(s) with respect to  the parameter using  bifurcation theory and degree theory.
	\end{abstract}
	\section{Introduction}
	We consider the following  elliptic equation with nonlinear boundary condition
	\begin{equation}
		\label{pde}
		\left\{
		\begin{array}{rcll}
			-\Delta u +u &=&  0 \quad &\mbox{in}\quad \Om\,;\\
			\frac{\p u}{\p \eta} &=& \la f(u)\quad &\mbox{on}\quad \p \Om,
		\end{array}
		\right. 
	\end{equation}
	where $\Om \subset \R^{N} (N \geq 2)$ is a bounded domain with $C^{2,\alpha}$ ($0< \alpha < 1$) boundary $\p \Om$,  $\p/\p\eta
	:=\eta(x)\cdot\nabla$ denotes the outer normal derivative on
	$\p\Om$,  $\la >0$ is a bifurcation parameter and  the nonlinearity on the boundary  
	$f :[0,\infty) \rightarrow [0,\infty) $ is locally Lipschtiz.
	\par Reaction-diffusion equations involving nonlinear boundary
	conditions, appear naturally in applications. For example, limb development which incorporates both outgrowth due to cell growth as well
	as cell division and interactions between morphogens produced in several very specific zones of the
	limb bud, see \cite{Dil_Oth_1999}. Another known application is when a highly exothermic reaction  takes place in a thin layer around a boundary \cite{Lacey-Ockendon-Sabina}, this information is then used in cryosurgery (surgery using local application of intense cold to destroy damaged tissue) \cite{Gurung-Gokul-Adhikary}. 
	\par Elliptic equations with nonlinear boundary conditions have been investigated extensively in recent years. Results on existence of positive solutions of problems with nonlinear boundary conditions can be found (without being exhaustive), using techniques such as,  monotone methods and functional analysis in \cite{Ama_1971, Inkmann}, concentration compactness method of Lions (see \cite{Lions-concentration-compactnessI, Lions-concentration-compactnessII}) in \cite{Garcia-Rossi-Sabina_2009}, bifurcation theory in \cite{Arrieta-Pardo-RBernal_2007, Arrieta-Pardo-RBernal_2009, Arrieta-Pardo-RBernal_2010, Liu-Shi_2018, Mavinga-Pardo_PRSE_2017,Pardo_Int_2012},  variational methods in \cite{Kajikiya-Naimen, Pap_Rad_2016, Quoirun-Umezu_2019}, and topological degree in \cite{Castro-Pardo_Inf_2017, Bon-Ros_2001}.

	\par Regarding  the nonlinear eigenproblem \eqref{pde}, we refer to \cite{Garcia-Peral-Rossi_2004, Quoirun-Umezu_2016, Pap_Rad_2016, Quoirun-Umezu_2019, Liu-Shi_2018} where there are  existence results with a parameter $\la$ on the boundary for a pure power sublinear nonlinearity. When $f$ is also a pure power and superlinear, we mention for instance \cite{Garcia-Peral-Rossi_2004} with a combination of interior and boundary reaction terms, and  also \cite{Hu-Yin_1994} where the authors describe the profile near blowup time for solution of the  associated parabolic problem. 
	\par 
	To the best of our knowledge, there are so far no existence results with respect to the parameter $\la$ on the boundary of problem \eqref{pde}  when the boundary nonlinearity $f$ is superlinear and subcritical, but not necessarily a pure power.  In this paper, we fill this gap by showing that there exists a positive weak solution for $\la$ small (see Theorem~\ref{thm:main}), depending only on the behavior of $f$ at infinity.
	Further, by imposing additional conditions on $f$ to guarantee bifurcation from the trivial solution, nonexistence for large $\la$,
	and some necessary technical assumptions,  we obtain  global bifurcation  and multiplicity results (see Theorem~\ref{th:conexion:0:inf}).
	\par 
	Main focus of this paper is to study \eqref{pde} when 
	$f$ is  superlinear and subcritical, that is,
	there exists a constant $b>0$ such that
	\begin{equation*}
		\hspace{-1cm} 
		\mathrm{(H)}_\infty \qquad\lim\limits_{s \to +\infty}\frac{f(s)}{s^p}=b
		\quad \mbox{ with }\quad \begin{cases}
			1 < p < \frac{N}{N-2} \; &\mbox{ if }\;N\ge 3,\\ 
			p>1\; &\mbox{ if }\;N=2   \,.
		\end{cases}
	\end{equation*}
	
	By a weak solution of problem \eqref{pde}, we mean a pair $(\la, u) \in (0, \infty) \times H^{1}(\Om)$ such that 
	\begin{equation}\label{weak:sol}
		\int_{\Om}\nabla u \nabla \psi + \int_{\Om}u\psi = \la \int_{\p \Om}f(u)\psi,\qquad \text{for all }\ \psi \in H^1(\Om). 
	\end{equation}
	Moreover,  
	one gets that the weak solution $u$ is  actually in $C^{2,\alpha}(\Om)\cap C^{1,\alpha}(\overline\Om)$ if $f$ satisfies the condition  $\mathrm{(H)}_\infty$ (see Corollary~\ref{cor:reg}).  Therefore, we use $C(\overline\Om)$ as our underlying space, and define the solution set as 
	$$
	\Sigma:= \left\{(\la, u) \in [0, +\infty) \times C(\overline\Om): (\la, u) \mbox{ is a weak solution of \eqref{pde}}\right\}\,.
	$$
	The closure of the set of nontrivial solutions will be denoted by $\mathscr{S}$.
	Let $\la_{0}, \la_{\infty}  \in [0, +\infty)$. We say that $(\la_{0}, 0) $ $ \big(\text{respectively},\ (\la_{\infty}, \infty)\big)$ is a bifurcation point from the  trivial solution  (respectively, from infinity) if there exists a sequence $(\la_n, u_n) \in \Sigma$ such that $\la_n \to \la_{0}$ and $\|u_n\|_{C(\overline{\Om})} \to 0$ (respectively, $\la_n \to \la_{\infty}$ and $\|u_n\|_{C(\overline{\Om})} \to +\infty$) as $n \to +\infty$.  
	Likewise, we say that a connected component  $\mathscr{C}$ bifurcates from the trivial solution at $(\la_{0}, 0)$ if $\mathscr{C}$ is a maximal closed connected subset of $
	\mathscr{S} 
	\cup (\la_{0}, 0)$ containing  $(\la_{0}, 0)$. 
	A connected component bifurcating from infinity can be defined similarly.
	\par We state our first result on local bifurcation from infinity.
	
	\begin{thm}
		\label{thm:main}
		Assume that $f$  satisfies $\mathrm{(H)}_\infty$. Then, there exists $\hat{\la} >0$ such that for all $\la \in (0, \hat{\la}]$, \eqref{pde} has a positive weak solution $u$ such that $\|u \|_{C(\overline{\Om})} \rightarrow \infty$ as $\la \rightarrow 0^+$. Moreover, there exists a connected component $\mathscr{C}^+ \subset \Sigma$, 
		of positive weak solutions of \eqref{pde}, bifurcating from infinity at $\la=0$, such that $\la$ takes all values in $(0, \hat{\la}]$ along $\mathscr{C}^+$ (see Fig.~\ref{fig:bif:infty}). 
	\end{thm}
	\setlength{\unitlength}{1cm}
	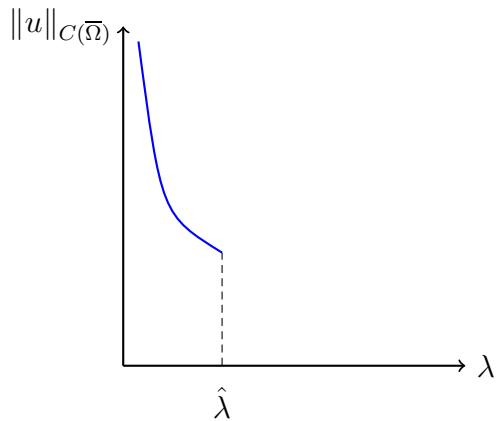
\begin{figure}[h]
		\centering
		\begin{tikzpicture}%
			\draw[thick,->] (0,0) -- (4.5,0) node[anchor=west] {$\lambda$};%
			\draw[thick,->] (0,0) -- (0,4.5) node[anchor=east] {$\|u\|_{C(\overline{\Omega})}$};%
			\draw[blue, thick] (0.2,4.3) .. controls (.5,2) .. (1.3,1.5);
			\draw[densely dashed] (1.3,0) -- (1.3, 1.5);
			\node[label=below:{{$\hat{\lambda}$}}] at (1.3,0) {};
		\end{tikzpicture}
		\caption{\small Bifurcation from infinity at $\lambda_{\infty}=0$}
		\label{fig:bif:infty}
	\end{figure}
	We  use uniform a-priori bounds results for asymptotically superlinear, subcritical 
	nonlinearities, and re-scaling argument together with degree theory and bifurcation theory to prove Theorem \ref{thm:main}. Such method was first used in \cite{Amb-Arc-Buf_1994} for a result in the Dirichlet boundary condition case.  We remark that Theorem~\ref{thm:main} is independent of the behavior of the nonlinearity $f$ away from infinity.
	
	\par Next, in order to discuss global bifurcation and multiplicity results, we
	impose additional conditions on the nonlinearity $f$. 
	First, we assume conditions on $f$ that guarantees bifurcation from the trivial solution, that is,  $f\in C^1([0,\infty))$ satisfies the following:
	\begin{equation*}
		\label{f(0):0:R}
		\hspace{-.6cm} 
		\mathrm{(H)}_0 
		\begin{cases}
			f(0)=0,\quad  f'(0)> 0,\\
			\text{and  there exists 
				a constant $\nu>1$  such that}
			\\
			f(s)=f'(0)s+\mathcal{R}(s) \text{ for}\ s\ge 0 \text{ with } \mathcal{R}(s)=O(s^\nu) \text{ as } s\to 0.
		\end{cases}
	\end{equation*}
	
	\par Second, to discuss the  bifurcation direction 
	of weak solutions near the bifurcation point,  following quantities play a crucial role.
	For $\nu>1$ as defined in $\mathrm{(H)}_0$, set
	\begin{equation}\label{R}
		\underline{\mathcal R}_0:=\liminf_{s\to 0^+} \frac{\mathcal{R}(s)}{s^\nu}
		\mbox{ and }
		\overline{\mathcal R}_0:=\limsup_{s\to 0^+} \frac{\mathcal{R}(s)}{s^\nu}
		\,. 
	\end{equation}
	
	\par Finally, let $\mu_1>0$ be the first Steklov eigenvalue and $ \varphi_1 \in H^1(\Om)$ the  corresponding nonnegative eigenfunction associated with the Steklov eigenvalue problem
	\begin{equation}\label{steklov}
		\left\{
		\begin{array}{rcll}
			-\Delta \psi +\psi &=&  0 \quad &\mbox{in}\quad \Om\,;\\
			\frac{\p \psi}{\p \eta} &=& \mu\psi\quad &\mbox{on}\quad \p \Om\,.
		\end{array}
		\right. 
	\end{equation}
	See Remark~\ref{rem:eig:fn} for the regularity and positivity of the nonnegative eigenfunction $\varphi_1$.
	
	\par Now, we state the following theorem concerning global bifurcation and multiplicity result.
	
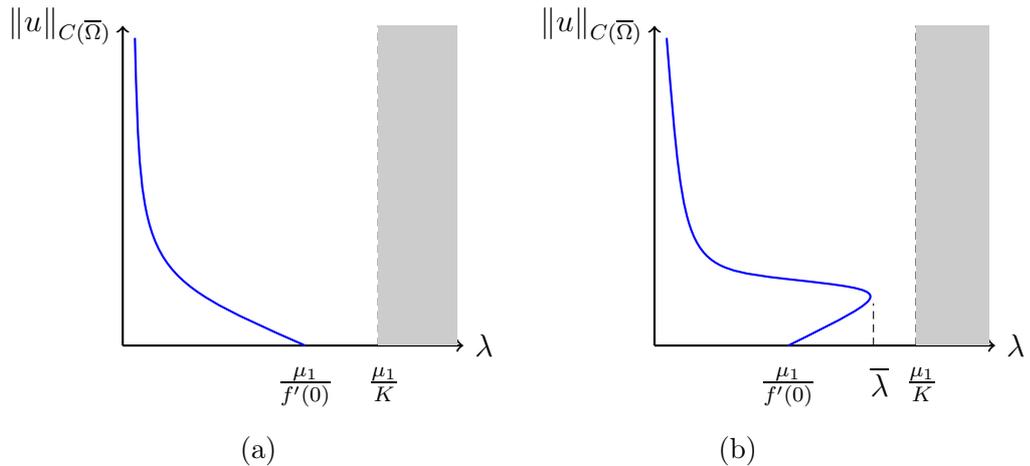
\begin{figure}[h]
		\setlength{\unitlength}{.8cm}
		\centering 
		\begin{subfigure}[b]{0.5\linewidth}
			\begin{tikzpicture}[scale=.8]
				\draw[thick,->] (0,0) -- (5.6,0) node[anchor=west] {$\lambda$};
				\draw[thick,->] (0,0) -- (0,5.3) node[anchor=east] {$\|u\|_{C(\overline{\Omega})}$}; 
				\draw[blue, thick] (0.2,5.1) .. controls (.3,1.2) .. (3, 0);
				\node[label=below:{$\frac{\mu_1}{f'(0)}$}] at (3, 0) {};
				\draw[densely dashed] (4.2,0) -- (4.2, 5.3);
				\node[label=below:{{$\frac{\mu_1}{K}$}}] at (4.3, 0) {};
				\fill[gray!40!white] (4.2, 0) rectangle (5.5,5.3);
			\end{tikzpicture}
			\caption{}  
		\end{subfigure}
		\begin{subfigure}[b]{0.4\linewidth}
			\begin{tikzpicture}[scale=.8]
				\draw[thick,->] (0,0) -- (5.6,0) node[anchor=west] {$\lambda$};
				\draw[thick,->] (0,0) -- (0,5.3) node[anchor=east] {$\|u\|_{C(\overline{\Omega})}$}; 
				\draw[blue, thick] (0.2,5.1) .. controls (.5,1.3) .. (2.2, 1.1);
				\draw[blue, thick] (2.2, 1.1) .. controls (4,.9) .. (2.2, 0);
				\node[label=below:{$\frac{\mu_1}{f'(0)}$}] at (2.2, 0) {};
				\draw[densely dashed] (3.6,0) -- (3.6,.7);
				\node[label=below:{{${\overline{\la}}$}}] at (3.7, 0) {};
				\draw[densely dashed] (4.3,0) -- (4.3, 5.3);
				\node[label=below:{{$\frac{\mu_1}{K}$}}] at (4.4, 0) {};
				\fill[gray!40!white] (4.3, 0) rectangle (5.5,5.3);
			\end{tikzpicture}
			\caption{}  
		\end{subfigure}
		\caption{\small Possible global bifurcation diagrams: (a) represents subcritical bifurcation; (b) represents supercritical bifurcation.  
		}
		\label{fig:sub:super:critical}
	\end{figure}

	\begin{thm}\label{th:conexion:0:inf}
		Let  $f\in C^1([0,\infty))$  be such that 
		hypotheses $\mathrm{(H)}_0$, $\mathrm{(H)}_\infty$  are satisfied. Suppose that there exists   $K>0$ such that 
		\begin{equation}
			\label{hyp:nonexist}
			f(s) \geq K s\quad \mbox{ for } s \geq 0\,. 
		\end{equation}
		\par Then, there exists a connected component  $ \mathscr{C}^+$ of positive weak solutions of \eqref{pde} emanating from the trivial solution at the bifurcation point
		$\big(\frac{\mu_1}{f'(0)},0\big)\in \Sigma$ possessing a unique bifurcation point from infinity at $\la=0$. More precisely, if $ (\la, u_{\la})\in \mathscr{C}^+$, then the following holds:
		\begin{equation}\label{unbdd:zero}
			\begin{cases}
				\|u_\la\|_{C(\overline{\Om})} \to 0 \ \ \qq{as} \la\to\frac{\mu_1}{f'(0)},\ \\
				\|u_\la\|_{C(\overline{\Om})} \to \infty \qq{as} \la\to  0^+\, \text{and} \\
				\text{if $(\la_\infty,\infty)$ is a bifurcation point from infinity, then $\la_\infty= 0$.}
			\end{cases}
		\end{equation}
		Moreover, problem \eqref{pde} has a  positive weak solution for any $\la \in\big(0, \frac{\mu_1}{f'(0)}\big)$ and no positive weak solutions for $\la > \frac{\mu_1}{K}$, see Fig.~\ref{fig:sub:super:critical}~(a)-(b).
		
		Furthermore, if  $\overline{\mathcal R}_0<0$, then the bifurcation from the trivial solution at  $\big(\frac{\mu_1}{f'(0)},0\big)$ is supercritical.
		In addition, there exists $\bar{\la}>\frac{\mu_1}{f'(0)}$ such that  problem \eqref{pde} 
		has at least two positive weak solutions for any $\la \in\big(\frac{\mu_1}{f'(0)},\bar{\la} \big)$, and at least one positive weak solution for  $\la =\frac{\mu_1}{f'(0)}$, and for  $\la =\bar{\la} $, as depicted in Fig.~\ref{fig:sub:super:critical} (b).
	\end{thm}
	
	\par We use novel approach of combining re-scaling argument used in the proof of local result, Theorem~\ref{thm:main}, and uniform a-priori bound together with the abstract local and global bifurcation theory \cite{C-R, R71}, and degree theory to prove Theorem~\ref{th:conexion:0:inf}.
	
	Analogous  existence results, such as  Theorems \ref{thm:main} and \ref{th:conexion:0:inf}, when the nonlinear reaction term appears as a source in the equation complemented with homogeneous Dirichlet
	boundary condition,  can be found, among others, in  \cite{Amb_Arc_book, Amb-Arc-Buf_1994} and the survey paper \cite{Lions}.

	\par Sections \ref{prelim} deals with some preliminaries such as regularity of weak solutions,  positivity,  and uniform a-priori bounds. In Section \ref{sec:proof:thm:main}, we prove Theorem~\ref{thm:main} using re-scaling argument and degree theory. 
	In Section \ref{sec:proof:thm:bif}, we collect results concerning bifurcation from the trivial solution using Rabinowitz's global bifurcation Theorem \cite{R71}. We also characterize the subcritical or supercritical nature of weak solutions near the bifurcation point.
	Finally, we prove Theorem~\ref{th:conexion:0:inf}  by combining bifurcation theory and degree theory.
	\par Unless otherwise specified,  solutions in this paper are understood as weak solutions, as defined in \eqref{weak:sol}. 
	\section{Preliminaries and Auxiliary Results}
	In this section, we discuss regularity and positivity of weak solutions of \eqref{pde}, and uniform a-priori bound result.
	\label{prelim}
	\subsection{Regularity of weak solutions and positivity}\label{regularity}
	Here, we state and prove regularity results for some linear and nonlinear problems,  which are relevant for our purposes. In particular, we prove that any weak solution of \eqref{pde} is in fact Hölder continuous, see Proposition  \ref{pro:reg} and  Corollary \ref{cor:reg}.
	
	To analyze the existence and regularity of weak solutions of  \eqref{pde}, we must set up the appropriate functional framework. To this end, we
	consider the following linear problem 
	\begin{equation}
		\label{lbp} 
		\left\{ \begin{array}{rcll}
			-\Delta  v +v &=&0 & \qquad \mbox{in } \Om\,;  \\
			\frac{\p   v}{\p \eta}&=&  \mathit{h} &  \qquad \mbox{on } \p
			\Om\,,
		\end{array}\right.
	\end{equation}
	where $h \in L^q(\p\Omega)$ for $q \geq 1$. It is known  that for each $q \geq 1$, \eqref{lbp} has a unique solution in $ W^{1,m}(\Om)$ and 
	\begin{equation}
		\label{ineq:m:h}
		\|v\|_{W^{1,m}(\Om)}\le  C\|\mathit{h}\|_{L^{q}(\p\Om)},\qquad \text{where }\quad 1 \leq  m\le Nq/(N-1)\,,
	\end{equation}
	see, for instance \cite{Mavinga-Pardo_PRSE_2017} for more details. We denote the solution operator corresponding to \eqref{lbp} by  
	\begin{equation*}
		\label{op:lin}
		T: L^{q}(\p\Om)\to W^{1,m}(\Om) \mbox{ with }  \mathit{T}\,\mathit{h}:=v\,.
	\end{equation*}
	\smallskip
	It is known that the trace operator
	\begin{equation}\label{op:trace}
		\Gamma : W^{1,m}(\Om ) \rightarrow L^{r}(\p
		\Om ) 
	\end{equation}
	is a continuous linear operator for every $r$ satisfying 
	$
	\frac{N-1}{r}\ge \frac{N}{m}-1$, and compact  if $\frac{N-1}{r}> \frac{N}{m}-1$, see \cite[Ch.~ 6]{KufnerJohnFucik77}.
	\smallskip
	Now, we define  the resolvent operator (also known as the {\it Neumann-to-Dirichlet operator}) $S:=T\circ \Gamma$ as
	\begin{equation}\label{def:S}
		S: L^q(\p\Om) \to L^r(\p\Om) \mbox{ given by }
		S\mathit{h}:=  \Gamma \big(\mathit{T}\,\mathit{h}\big)=\Gamma v\,,
	\end{equation}
	for any $q \geq 1$ and for all $r$ satisfying $\frac{N-1}{r}\ge \frac{N-m}{m}$ with $1 \leq m\le Nq/(N-1)$, given schematically by
	$$
	L^q(\p\Om)\stackrel{T}{\longrightarrow} W^{1,m}(\Om)\stackrel{\Gamma}{\longrightarrow}L^r(\p\Om)\,.
	$$
	Note that if $\frac{N-1}{r}> \frac{N-m}{m}$ then $S$ is compact by the compactness of $\Gamma$.
	\smallskip
	The following Lemma states the regularity of the solution of the linear problem \eqref{lbp}. In particular, if  ${q}>N-1$, then  $ v \in C^\alpha(\overline\Om)$.
	\begin{lem}\label{lem:reg:lbp}
		Let $N\ge  2$ and $\mathit{h}\in L^q (\p \Om )$
		with $q\ge  1$. Then, the unique solution $v=\mathit{T}\mathit{h}$ of the linear problem \eqref{lbp} satisfies the following:
		\begin{enumerate}
			\item [(i)] If $1\le  q<N-1$, then $\Gamma v\in L^{r}(\p\Om)$
			for all $1\le  r\le  \frac{(N-1)q}{N-1-q}$ and the map $
			S:L^q(\p\Om)\to L^r(\p\Om)$ is
			continuous for $1\le  r\le  \frac{(N-1)q}{N-1-q}$  and compact for 
			$1\le  r<\frac{(N-1)q}{N-1-q}$. 
			
			\item [(ii)] If $q=N-1$, then $\Gamma v\in L^{r}(\p\Om)$ for 
			all $r\ge  1$ 
			and the map $ S:L^q(\p\Om)\to L^r(\p\Om)$ is
			continuous and compact for $1\le  r<\infty$. 
			
			\item [(iii)] If $q>N-1$, then $v\in C^{\alpha }(\overline{\Om})$ with
			$\|v\|_{C^{\alpha }(\overline{\Om})}\le 
			C\|h\|_{L^q(\p\Om)}$ for 
			some $\alpha\in (0,1)$. Moreover, 
			$\Gamma v\in C^{\alpha}(\p\Om)$ 
			and the map $S:L^q(\p\Om)\to C^\alpha(\p\Om)$ is  continuous and compact.
			
			\item [(iv)] If $\mathit{h}\in C^{\alpha }(\p\Om)$, then $v\in C^{2,\alpha }(\Om)\cap C^{1,\alpha }(\overline{\Om})$.
		\end{enumerate}
	\end{lem}
	\begin{proof}
		See \cite[Lemma 2.1]{Arrieta-Pardo-RBernal_2007} for proofs of items (i)-(iii).
		\par (iv) Since $h \in C^{\alpha}(\p \Om)$ then  by (iii) $v \in C^\alpha(\overline{\Om})$,  
		it follows from the first equation  in  \eqref{lbp} that $v \in C^{2, \alpha}(\Om)$. Furthermore, using the second equation of \eqref{lbp}, it follows that $v\in C^{2,\alpha}(\Om)\cap C^{1,\alpha}(\overline\Om)$.
		
	\end{proof}
	\par In what follows, we will show that any weak solution $u$ of our nonlinear problem \eqref{pde} lies in fact in $C^\alpha(\overline\Om)$  for some  $\alpha \in (0, 1)$. To
	accomplish this, we will establish regularity results for problems with nonlinearities satisfying $\mathrm{(H)}_\infty$.
	Hereafter, we will use the same symbol to denote both the function and the associated Nemytskii operator.
	\begin{pro}\label{pro:reg}
		Let $N\ge  2$ and $h :[0,\infty) \rightarrow [0,\infty) $ be locally Lipschitz continuous satisfying condition $\mathrm{(H)}_\infty$. Let
		$v$ be a nontrivial weak solution of the following problem 	
		\begin{equation*}
			\label{pde2}
			\left\{
			\begin{array}{rcll}
				-\Delta v +v &=&  0 \quad &\mbox{in}\quad \Om\,;\\
				\frac{\p v}{\p \eta} &=&  h(v)\quad &\mbox{on}\quad \p \Om\,.
			\end{array}
			\right. 
		\end{equation*}	
		Then, 
		\begin{equation*}\label{regestimate:v:g}
			||v||_{C^\alpha(\overline{\Om})}\le C(1+||\Gamma v||_{L^{r}(\p\Om)})
		\end{equation*}
		for some positive $\alpha \in (0, 1)$, where  $\displaystyle r=\frac{2(N-1)}{N-2} $ if $N>2$, and $r\ge 1$ when $N=2$.
	\end{pro}
	
	\begin{proof}
		We assume $N>2$, since  the proof is trivial when $N=2$.  By definition of a weak solution and the trace operator, \eqref{op:trace}, $v\in H^1(\Om)$ and its  trace $\Gamma v\in L^{r}(\p \Om )$, where $1 \leq r\le r_0:=\frac{2(N-1)}{N-2}$, respectively. 
		It follows from the  condition $\mathrm{(H)}_\infty$ that
		\begin{equation}
			\label{h:subcrit}
			h\big(\Gamma v\big)\le C(1+|\Gamma v|^p),
		\end{equation}
		and by the continuity of the Nemytskii operator
		$$
		h\big(\Gamma v\big)\in  L^{q_0}(\p\Om),\qquad\text{where}\qquad  q_0:=\frac{r_0}{p}=\frac{2(N-1)}{p(N-2)}.
		$$
		Now we proceed with the  bootstrap argument.  For  $h\big(\Gamma v\big)\in  L^{q_0}(\p\Om)$ and \eqref{ineq:m:h}, we have
		$$
		v\in W^{1,s_{i}}(\Om), \mbox { where }s_{i}:=\frac{N{q}_{i-1}}{N-1}\quad \mbox{ for }i=1,2,\ldots\,.$$ 
		By \eqref{op:trace}, we get 
		$$\Gamma v\in  L^{r_i}(\p\Om),  \mbox { where }  r_i
		:=\frac{(N-1)q_{i-1}}{N-1-q_{i-1}}\quad \mbox{ for }i=1,2,\ldots\,. $$
		Then, using \eqref{h:subcrit} and the continuity of the Nemytskii operator  $$ h\big(\Gamma v\big)\in  L^{q_i}(\p\Om),
		\mbox{ where } q_i:=\frac{r_i}{p}\quad \mbox{ for }i=1,2,\ldots\,.$$
		\par If $q_i  > N-1$ for some $i\in \N $, then $v\in C^{\alpha} (\overline{\Om})$  for some $\alpha\in (0,1)$ by Lemma \ref{lem:reg:lbp} (iii). 
		\par If $q_i = N-1$ for some $i\in \N $, then by Lemma~\ref{lem:reg:lbp} (ii), $\Gamma v\in L^r(\p\Om)$ for $r\ge1$.
		By \eqref{h:subcrit}, $h(\Gamma v)\in L^m$ for $m\ge 1$ . 
		Using the $L^q$-estimates for second-order linear elliptic equations, we get that $u$ is actually in $W^{1,s}(\Om)$ for any $s>1$, in particular for $s>N$. By the continuity of the embedding $W^{1,s}(\Om)\hookrightarrow C^{\alpha}(\overline{\Om})$ for $s>N$, one has that $v\in C^{\alpha}(\overline{\Om})$, see e.g \cite[p. 285]{Brezis_2011}. 
		
		Now suppose $q_i<N-1$. 
		Then, 
		$$
		\frac{1}{r_1}= \frac{1}{q_0} - \frac{1}{N-1}
		=\frac{p(N-2)-2}{2(N-1)}<\frac{N-2}{2(N-1)}=\frac{1}{r_0}\qquad\text{iff}\ 
		p<\frac{N}{N-2}\,.
		$$
		If  $r_{i}>r_{i-1}$, then 
		$$
		\frac{1}{r_{i+1}}= \frac{1}{q_i} - \frac{1}{N-1}
		= \frac{p}{r_i} - \frac{1}{N-1}< \frac{p}{r_{i-1}} - \frac{1}{N-1}
		=\frac{1}{q_{i-1}} - \frac{1}{N-1}=\frac{1}{r_{i}}\,.
		$$
		Hence, by induction $\{r_i\}$ is strictly increasing. Then, clearly $\{s_i\}$ and $\{q_i\}$ are strictly increasing as well.
		\smallskip
		
		Suppose  $q_i<N-1$ for all $i\in\mathbb{N}$. Since 
		$\{ q_i \}$ is  strictly increasing  and $1\le q_i<N-1$ for all $i\in\mathbb{N}$,  $\displaystyle q_i \to q_\infty $ for some  $1\le q_\infty \le N-1$. If $ q_\infty = N-1$, then fixing $\varepsilon>0$, there exists an $i_{0}\in\N$ such that $N-1>q_{i_0}\ge N-1-\varepsilon$. However, $\{q_i\}$ is strictly increasing, hence $q_{i_{0}+1}>q_{i_0} \ge N-1$,  a contradiction. As a consequence, $ q_\infty < N-1$.
		Define $r_\infty:= \displaystyle \lim_{i\to \infty}r_i = \displaystyle \lim_{i\to \infty}pq_i = p q_{\infty} >0 $. Note that 
		\begin{align*}
			q_{i+1}-q_i&=\frac{r_{i+1}-r_i}{p}
			=\frac{r_{i+1}\,r_i}{p}\left(\frac{1}{r_{i}}-\frac{1}{r_{i+1}}\right)\\
			&=\frac{r_{i+1}\,r_i}{p}\left(\frac{1}{q_{i-1}}-\frac{1}{q_{i}}\right)=\frac{r_{i+1}\,r_i}{p}\, \frac{q_{i}-q_{i-1}}{q_{i-1}\,q_{i}}
		\end{align*}
		hence
		\begin{align*}
			\frac{q_{i+1}-q_i}{q_{i}-q_{i-1}}&=\frac{r_{i+1}\,r_i}{p}\, \frac{1}{q_{i-1}\,q_{i}}=\frac{r_{i+1}\,r_i}{p}\, \frac{p^2}{r_{i-1}\,r_{i}}=p\,\frac{r_{i+1}}{r_{i-1}}\,.
		\end{align*}
		Taking the limit as $i$ goes to infinity and noting that $r_{\infty}>0$, we have
		\begin{align*}
			\lim_{i\to \infty} \frac{q_{i+1}-q_i}{q_{i}-q_{i-1}} &=p>1\,.
		\end{align*}
		This contradicts the boundedness of $\{q_i\}$.
		Therefore, there exists $i_0 \in \N$ such that $q_{i_0}\ge N-1$ and hence $v \in C^{\alpha}(\overline{\Om})$ for some $\alpha \in (0,1)$, as desired.
		Furthermore, the estimate in Lemma~\ref{lem:reg:lbp} and \eqref{h:subcrit} give  \begin{equation}\label{regestimate:v:h}
			||v||_{C^\alpha(\overline{\Om})}\le \|h(\Gamma v)\|_{L^{q_0}(\p\Om)}\le \|C(1+|\Gamma v|^p)\|_{L^{q_0}(\p\Om)} \le C(1+||\Gamma v||_{L^{r}(\p\Om)})
		\end{equation}
		where  $\displaystyle r=pq_0=\frac{2(N-1)}{N-2} $ if $N>2$.
	\end{proof}

	\begin{cor}\label{cor:reg}
		\rm
		Assume that  the nonlinearity $f:[0,\infty) \rightarrow [0,\infty)$ is  locally Lipschitz continuous and satisfies condition $\mathrm{(H)}_\infty$. 
		Fix any $\La >0$ and let $u$ be a weak solution of the nonlinear problem \eqref{pde} for some $0< \la\le  \La $. 
		Then 
		\begin{equation*}\label{regestimate}
			||u||_{C^\alpha(\overline{\Om})}\le C(1+||\Gamma u||_{L^{r}(\p\Om)}),
		\end{equation*}
		for some $\alpha \in (0, 1)$ and some $C=C(\La)>0$, where  $\displaystyle r=\frac{2(N-1)}{N-2}$ for $N>2$ and $r\ge 1$ for $N=2$.	Moreover, $u\in C^{2,\alpha}(\Om)\cap C^{1,\alpha}(\overline\Om) $.
	\end{cor}
	
	\begin{proof}
		Proposition~\ref{pro:reg} yields the proof for the first part.

		\par Since  $u\in C^{\alpha}(\overline{\Om})$, $f$ is locally Lipschtiz continuous, $f(u) \in C^{\alpha}(\p \Om)$. The conclusion follows from  Lemma \ref{lem:reg:lbp} (iv).
	\end{proof}
	
	Under additional assumption on the nonlinearity $f$, Corollary \ref{cor:reg} can be rewritten in the following way.
	\begin{pro}\label{pro:0:inf}
		Assume that  the nonlinearity  $f\in C^1([0,\infty))$ satisfies conditions $\mathrm{(H)}_0$ and $\mathrm{(H)}_\infty$. For any fixed $\La >0$, if $u$ is  a  weak solution of the nonlinear problem \eqref{pde} for some $0\leq \la \le  \La $,
		then 
		\begin{equation*}\label{regestimate:0:inf}
			\|u\|_{C^\al (\overline{\Om})}\le C||\Gamma u||_{L^{r}(\p \Om)},
		\end{equation*}
		for some $\al \in (0, 1)$, where $C=C(\La )$ and $ r=\frac{2(N-1)}{N-2} $  if $N>2$, and $r\ge 1$ when $N=2$.
	\end{pro}
	
	\begin{proof}
		
		\par Note that under  conditions $\mathrm{(H)}_0$ and $\mathrm{(H)}_\infty$, for any $\e>0$, there exists a constant $C_\e>0$ such that
		\begin{equation*}\label{f:0:inf:eps}
			f(s)\le (1+\e)f'(0)|s|+C_\e|s|^p.
		\end{equation*}
		In particular, there exists a constant $C>0$ such that
		$
		f(s)\le C(|s|+|s|^p). 
		$ 
		Hence, the conclusion follows from \eqref{regestimate:v:h}.
	\end{proof}
	Next lemma shows that any nonnegative nontrivial solution of \eqref{lbp} is positive on $\overline{\Om}$.
	\begin{pro}\label{pro:v>0}
		Let  $v \in C^2(\Omega) \cap C^1(\overline{\Omega})$ be a solution of \eqref{lbp} for $h \geq 0$ with $h \not\equiv 0$. Then $v > 0$ on $\overline{\Omega}$.
	\end{pro}

	\begin{proof}
		Clearly $v>0$ in $\Om$ by the strong Maximum Principle,  see \cite[p.~127]{Ama_1971}.
		Assume to the contrary  that there exists an $x_0\in\p\Om$ such that $v(x_0)=0$. 	By the Hopf's Lemma 
		(\cite[Lem.~3.4]{Gil-Trud_2001})
		$
		\frac{\p v}{\p \eta}(x_0)<0,
		$ 
		contradicting the boundary condition $\frac{\p v}{\p \eta}(x_0)=h(x_0) \geq 0$. As a conclusion, $v>0$ for all $ x\in\overline{\Om}$.
	\end{proof}
	
	\begin{rem}
		{\rm 
			Let $f:[0, \infty) \to [0, \infty)$ be a locally Lipschitz continuous satisfying condition $\mathrm{(H)}_\infty$ and $u$ be a  weak solution of \eqref{pde} for some $\lambda>0$. Then,  Corollary~\ref{cor:reg} implies that $u\in C^{2,\alpha}(\Om)\cap C^{1,\alpha}(\overline\Om)$ and hence  $u>0$ on $\overline{\Omega}$ by Proposition~\ref{pro:v>0}. }
	\end{rem}
	\begin{rem}\label{rem:trace:u:reg}{\rm
			Let $f:[0, \infty) \to [0, \infty)$ be a locally Lipschitz continuous function satisfying condition $\mathrm{(H)}_\infty$. Then, for a given $u\in C(\overline\Om)$,   
			$f(\Gamma u)$ maps $C(\overline\Om)$ into  $L^q(\p\Om)$ with $q>1$ by the continuity of the Nemytskii operator associated with $f$, see \cite[Lemma 3.1]{Ama_1976-b}. Then, using \eqref{def:S}, we have that 
			\begin{equation}\label{mapping:S:f} 
				S\circ f\circ \Gamma: C(\overline\Om)\longrightarrow L^r(\p\Om)\stackrel{{\rm Cor}. \ref{cor:reg}}{ \longrightarrow } C^\alpha(\overline\Om)\stackrel{c}{\hookrightarrow} C(\overline\Om)\,,\end{equation}
			is compact, and $v=(S\circ f\circ \Gamma) u$ is the weak solution of 
			\begin{equation*}
				\left\{
				\begin{array}{rcll}
					-\Delta v +v &=&  0 \quad &\mbox{in}\quad \Om\,;\\
					\frac{\p v}{\p \eta} &=& \la f(u)\quad &\mbox{on}\quad \p \Om\,.
				\end{array}
				\right. 
			\end{equation*}
			More precisely,  
			\begin{equation*}\label{fixed:point}
				u\text{ is  a weak solution of   \eqref{pde} for } \la >0 \iff
				u =\la   S\big( f(\Gamma u)\big)\,.
			\end{equation*} 
		}
	\end{rem}
	
	\smallskip
	We end this subsection with a remark about the sign and regularity of the eigenfunction $\varphi_1$ corresponding to the first Steklov eigenvalue $\mu_1$ of problem \eqref{steklov}. 
	\begin{rem}
		\label{rem:eig:fn}
		{\rm 
			By the regularity of weak solutions, see \cite[Thm.~9.26]{Brezis_2011},
			and repeating the arguments as in the proof of Corollary \ref{cor:reg}, the eigenfunction $\varphi_1$ corresponding to the first Steklov eigenvalue $\mu_1$ of problem \eqref{steklov} is  in  $ C^{2, \alpha}(\Om)\cap C^{1, \alpha}(\overline\Om)$ with $0<\alpha<1$ (see also e.g.   \cite[Thm~8.12]{Mav_2012}).  
			Therefore, by Proposition~\ref{pro:v>0}   $\varphi_1>0$ on $\overline\Om$.
		}
	\end{rem}
	
	\subsection{Uniform a-priori bound}
	\par Our main tool in the proof of Theorem~\ref{thm:main} is degree theory, for which the following uniform a-priori bound is crucial.
	To state the result, consider
	\begin{equation}
		\label{pde:apriori}
		\left\{
		\begin{array}{rcll}
			-\Delta u +u &=&  0 \quad &\mbox{in}\quad \Om\,;\\
			\frac{\p u}{\p \eta} &=& b u^p + \z(x, u)\quad &\mbox{on}\quad \p \Om\,,
		\end{array}
		\right. 
	\end{equation}
	where $p$ is as in $\mathrm{(H)}_\infty$, and for a.e. $x \in \overline{\Om}$, all $\sigma \in \R$,   and 
	
	\begin{equation}\label{o}
		\lim_{\s\to\infty }\frac{|\z(x, \sigma)|}{|\sigma|^p}\, =0.
	\end{equation}
	
	\par While we are not aware of any paper that establishes uniform a-priori  estimate for \eqref{pde:apriori}, the result below 
	follows by adapting the proof for systems case in \cite[Thm.~3.7]{Bon-Ros_2001}. Their proof is written for 
	$|\z(x, \sigma)| \leq c(1+ |\sigma|^r) \, \mbox{ for some } 0<r<p$, but the same arguments can be used to prove the existence of a priori bound under condition \eqref{o}.
	\begin{pro}
		\label{prop:apriori}
		There exists a constant $M > 0$ such that every positive  solution $u\in C(\overline{\Om})$ of \eqref{pde:apriori} satisfies
		\begin{equation*}
			\label{eq:apriori}
			\| u \|_{C(\overline{\Om})} \leq M\,.
		\end{equation*}
	\end{pro}
	
	\smallskip
	\section{Proof of Theorem~\ref{thm:main}}
	\label{sec:proof:thm:main}
	Our proof is motivated by \cite{Amb-Arc-Buf_1994}. In particular, we re-scale \eqref{pde} in such a way that the transformed problem approaches a limiting problem of "pure power type" as $\la \to 0^+$. Then, using $\la \geq 0$ as the homotopy parameter, we obtain a positive weak solution of the re-scaled problem, hence of \eqref{thm:main} for $\la> 0$ small.

	\par First, let us  extend $f$ to $\R$ by setting $f(s)=f(|s|)$ for $s \in \R$. Now consider the problem
	\begin{equation}
		\label{pde:extend}
		\left\{
		\begin{array}{rcll}
			-\Delta u +u &=&  0 \quad &\mbox{in}\quad \Om\,;\\
			\frac{\p u}{\p \eta} &=& \la f(|u|) \quad &\mbox{on}\quad \p \Om\,.
		\end{array}
		\right. 
	\end{equation}
	Note that for $\la >0$, $u$ is a solution of \eqref{pde:extend} if and only if $w=\la^{\frac{1}{p-1}}u$ satisfies
	\begin{equation}
		\label{pde:extend:2}
		\left\{
		\begin{array}{rcll}
			-\Delta w +w &=&  0 \quad &\mbox{in}\quad \Om\,;\\
			\frac{\p w}{\p \eta} &=& \la^{\frac{p}{p-1}} f\big(\la^{-\frac{1}{p-1}}|w|\big) \quad &\mbox{on}\quad \p \Om\,.
		\end{array}
		\right. 
	\end{equation}
	For $\la >0$, define 
	\begin{align*}\label{tilde:f}
		\tilde{f}(\la, s) 
		&:=\la^{\frac{p}{p-1}} f\big(\la^{-\frac{1}{p-1}}|s|\big)\\
		&= \la^{\frac{p}{p-1}} \left[f\big(\la^{-\frac{1}{p-1}}|s|\big) 
		- b\big(\la^{-\frac{1}{p-1}}|s|\big)^p\right] + b |s|^p\,. \nonumber
	\end{align*}
	We observe that 
	$$
	\lim\limits_{\substack{\la \to 0^+\\ s \to s_0}} \la^{\frac{p}{p-1}} f\big(\la^{-\frac{1}{p-1}}|s|\big) = b|s_0|^p\,,
	$$
	due to superlinear condition at infinity  $\mathrm{(H)}_\infty$ for $s_0 \neq 0$, and  by the continuity of $f$  at $s_0=0$. Therefore, we can define $\tilde{f}$ at $\la=0$ by setting $\tilde{f}(0, s):=b|s|^p$. Therefore, since $f$ is  Lipschitz continuous, so is 
	$\tilde{f}: [0, +\infty) \times \R \to [0, +\infty)$ defined above.

	\par Then the goal is to study the following re-scaled  problem for $\la \geq 0$
	\begin{equation}
		\label{pde:rescaled}
		\left\{
		\begin{array}{rcll}
			-\Delta w +w &=&  0 \quad &\mbox{in}\quad \Om\,;\\
			\frac{\p w}{\p \eta} &=&  \tilde{f}(\la, w) \quad &\mbox{on}\quad \p \Om\,,
		\end{array}
		\right. 
	\end{equation}
	while keeping in mind that \eqref{pde:rescaled} reduces to the limiting problem for $\la=0$
	\begin{equation}
		\label{pde:limit}
		\left\{
		\begin{array}{rcll}
			-\Delta w +w &=&  0 \quad &\mbox{in}\quad \Om\,;\\
			\frac{\p w}{\p \eta} &=& b |w|^p \quad &\mbox{on}\quad \p \Om\,.
		\end{array}
		\right. 
	\end{equation}
	
	\par Our strategy to proceed with the proof of Theorem~\ref{thm:main} is as follows: 1) we show that the limiting problem \eqref{pde:limit}, corresponding to $\la=0$, has a positive solution using the Leray-Schauder degree, 2) show that the  re-scaled problem \eqref{pde:rescaled} has a positive solution using  1) and  $\la \geq 0$ as the homotopy parameter, then 3) return to the original problem via the re-scaling.
	\par To set up for the Leray-Schauder degree, we formulate the problem \eqref{pde:rescaled} in an abstract setting in terms  of the  compact and Nemytskii  operators. 
	For this, we define the compact map $ \tilde{\mathcal{F}}: [0, +\infty) \times C(\overline{\Om}) \to  C(\overline{\Om})$ given by
	\begin{equation*}
		\tilde{\mathcal{F}}(\la, v):=S (\tilde{f}(\la, \Gamma (v)))\,,
	\end{equation*}
	where $\tilde{f}(\la, \cdot)$ denotes the Nemytskii operator corresponding to $\tilde{f}(\la, \cdot)$, and $S$ is as defined in Remark~\ref{rem:trace:u:reg}. It follows from Remark~\ref{rem:trace:u:reg}  that 
	$$
	(\la, w)\ \  \text{is a weak solution of  \eqref{pde:rescaled}} \iff \tilde{{\mathcal F}}(\la, w)=  w\,.  
	$$
	
	\par First we establish the following result regarding the limiting problem \eqref{pde:limit}.
	\begin{lem}
		\label{lem:0:r}
		There exists $r>0$ such that for all $\theta \in [0, 1]$ and all $w \in C(\overline \Om)$ with $\|w\|_{C(\overline\Om)}=r$, $w \neq \theta \tilde{{\mathcal F}}(0, w)$. Consequently $\text{deg}(I-\tilde{{\mathcal F}}(0, \cdot), B_r(0), 0)=1$.
	\end{lem}
	
	\begin{proof} 
		Suppose to the contrary that for each $r>0$, there exists $\theta \in[0, 1]$ such that the operator equation
		$$
		w =\theta \tilde{{\mathcal F}}(0, w) \quad\text{for}\quad \theta \in [0, 1]$$
		has a solution $w \in C(\overline \Om)$ with $\| w\|_{C(\overline\Om)}=r$, that is, $w$  is a solution of 
		\begin{equation}
			\label{pde:theta}
			\left\{
			\begin{array}{rcll}
				-\Delta w +w &=&  0 \quad &\mbox{in}\quad \Om\,;\\
				\frac{\p w}{\p \eta} &=& \theta\, b |w|^p\quad &\mbox{on}\quad \p \Om\,.
			\end{array}
			\right. 
		\end{equation}
		Clearly  $w\ne 0$  since $\| w\|_{C(\overline\Om)}=r>0$. Hence $w>0$
		in $\overline \Om$ by Proposition~\ref{pro:v>0}.
		
		\par Now, let $0< \varepsilon< \mu_1$ be fixed. Since $p>1$, there exists $r^*>0$ such that $bs^p <  \varepsilon s$ for $0 < s \leq  r^*$.
		Then there exists $\theta_{r^*} \in [0, 1]$ and     a solution $w_{r^*} > 0$  of \eqref{pde:theta} such that $\| w_{r^*}\|_{C(\overline\Om)}={r^*}$, and  
		$w_{r^*}$  satisfies
		$bw_{r^*}^p <  \varepsilon w_{r^*}$ whenever $\| w_{r^*}\|_{C(\overline\Om)}={r^*}$. Using $\varphi_1\geq 0$ as the test function and the fact that $\theta_{r^*} \in [0, 1]$, we have
		\begin{align*}
			0\ & =   \int_{\Om}\nabla w_{r^*} \nabla \varphi_1 + \int_{\Om}w_{r^*}\varphi_1 - \theta_{r^*} b \int_{\p \Om} w^p_{r^*}\varphi_1\\
			&\geq \int_{\Om}\nabla w_{r^*} \nabla \varphi_1 + \int_{\Om}w_{r^*}\varphi_1 - \varepsilon \int_{\p \Om} w_{r^*} \varphi_1
			= (\mu_1 -  \varepsilon)\int_{\p \Om} w_{r^*}\varphi_1\,,   
		\end{align*}
		a contradiction since $\varepsilon<\mu_1$.   Thus there exists $r>0$ such that for all $\theta \in [0, 1]$ and all $ w \in C(\overline \Om)$ with $\| w\|_{C(\overline\Om)}=r$, $w \neq \theta \tilde{{\mathcal F}}(0, w)$. Therefore, using $\theta \in [0, 1]$ as a homotopy parameter, we get
		$$
		\text{deg}(I-\tilde{{\mathcal F}}(0, \cdot), B_r(0), 0)=
		\text{deg}(I- \theta \tilde{{\mathcal F}}(0, \cdot), B_r(0), 0)
		=
		\text{deg}(I, B_r(0), 0)
		=1\,,
		$$
		as desired. This completes the proof of Lemma~\ref{lem:0:r}.
	\end{proof}
	
	\begin{lem}
		\label{lem:0:R}
		There exists $R > r>0$ and $0 \leq z \in C(\overline\Om)$ such that 
		$
		w \neq \tilde{{\mathcal F}}(0, w) + tz
		$
		for all $t \geq 0$ and all $w \in C(\overline\Om)$ with $\|w\|_{C(\overline\Om)}=R$. Consequently, $\text{deg}(I-\tilde{{\mathcal F}}(0, \cdot), B_R(0), 0)=0$.
	\end{lem}
	\begin{proof} Let 
		$0 \leq z \in C(\overline \Om)$ be the unique solution of 
		\begin{equation*}
			\label{pde:1}
			\left\{
			\begin{array}{rcll}
				-\Delta z +z &=&  0 \quad &\mbox{in}\quad \Om\,;\\
				\frac{\p z}{\p \eta} &=& 1\quad &\mbox{on}\quad \p \Om\,.
			\end{array}
			\right. 
		\end{equation*}
		Then, we observe that the operator equation 
		\begin{equation*}
			\label{op:eq:limit:t}
			w = \tilde{{\mathcal F}}(0,  w) + t\,z
		\end{equation*}
		corresponds to the PDE
		\begin{equation}
			\label{pde:limit:t}
			\left\{
			\begin{array}{rcll}
				-\Delta w +w &=&  0 \quad &\mbox{in}\quad \Om\,;\\
				\frac{\p w}{\p \eta} &=& b |w|^p +t\quad &\mbox{on}\quad \p \Om\,.
			\end{array}
			\right. 
		\end{equation}
		\noindent{\bf Step~1:} We show that there exists $t_0>0$ such that  \eqref{pde:limit:t} does not have a solution for $t \ge t_0$. 
		\par For this, let $\mu>\mu_1$ be fixed. Then there exists $t_0>0$ such that $bs^p + t > \mu s  + t- t_0$ for $t \geq 0$. Suppose by contradiction that there exists $t_1 \geq t_0$ such that $w \ge 0$ is a solution of \eqref{pde:limit:t}.  Using $\varphi_1\ge 0$ as the test function, we get 
		\begin{align*}
			0 &= \int_{\Om}\nabla w \nabla \varphi_1 + \int_{\Om}w\varphi_1 -  \int_{\p \Om}[b w^p + t_1]\varphi_1\\
			&= \int_{\p \Om}[b w^p + t_1]\varphi_1 - \int_{\Om}\nabla w \nabla \varphi_1 - \int_{\Om}w\varphi_1\\
			&> \int_{\p \Om}[\mu w + (t_1- t_0)]\varphi_1 - \int_{\Om}\nabla w \nabla \varphi_1 - \int_{\Om}w\varphi_1   \\
			&\geq   \mu\int_{\p \Om} w \varphi_1 - \int_{\Om}\nabla w \nabla \varphi_1 - \int_{\Om}w\varphi_1  
			= (\mu - \mu_1)\int_{\p \Om} w\varphi_1\,,
		\end{align*}
		which is a contradiction since $\mu> \mu_1$. This establishes Step~1, which implies that for all $a>0$, $w \neq \tilde{{\mathcal F}}(0, w) + t_0\,z$ for all $w \in C({\overline \Om})$ with $\|w\|_{C(\overline\Om)}=a$ for any $a>0$. Hence, for any $a>0$, we have
		\begin{equation}
			\label{eq:deg:t_0}
			\text{deg}(I-\tilde{{\mathcal F}}(0, w) + t_0\,z, B_a, 0)=0\,.
		\end{equation} 
		
		\noindent{\bf Step~2:}  We show there exists $R>r>0$ such that for all $t \in [0, t_0]$,
		$$\text{deg}(I-\tilde{{\mathcal F}}(0, w) + t\,z, B_R, 0)=0\,.
		$$
		Indeed, by Proposition~\ref{prop:apriori} with $\xi(t)\equiv t \in [0, t_0]$, there exists $M>0$ such that $  \|w\|_{C(\overline{\Om})} \leq  M$. By taking $R > \max\{r, M\}$, we get 
		$w \neq \tilde{{\mathcal F}}(0, w) + t\,z$ for all $w \in C({\overline \Om})$ with $\|w\|_{C(\overline\Om)}=R$ and $t \in [0, t_0]$. Then, using \eqref{eq:deg:t_0}, we get
		$$\text{deg}(I-\tilde{{\mathcal F}}(0, w), B_R, 0)=
		\text{deg}(I-\tilde{{\mathcal F}}(0, w) + t_0\,z, B_R, 0)
		=0\,,
		$$
		as desired, establishing Step~2. This completes the proof of Lemma~\ref{lem:0:R}.
	\end{proof}
	\par 
	Now we show that the limiting problem \eqref{pde:limit} has a positive solution.\\
	Indeed, it follows from Lemma~\ref{lem:0:r}, Lemma~\ref{lem:0:R} and the excision property of degree that 
	\begin{equation}
		\label{deg:0:-1}    
		\text{deg}(I-\tilde{{\mathcal F}}(0, w), B_R \setminus \overline{B}_r, 0)=-1\,.
	\end{equation}
	Therefore, there exists a solution of $w=\tilde{{\mathcal F}}(0, w)$, or equivalently a weak solution of  \eqref{pde:limit}, say $w_0 \in B_R \setminus \overline{B}_r$.
	Using the fact that $\|w_0\|_{C(\overline\Om)}>r >0$, it follows from Proposition~\ref{pro:v>0} that $w_0 > 0$ in $\overline{\Om}$.
	\medskip
	\par Now we use $\la\geq 0$ as homotopy parameter to establish the following existence result for the re-scaled problem \eqref{pde:rescaled}.
	\begin{lem}
		\label{lem:rescaled}
		There exists $\hat{\la}>0$ such that 
		\begin{description}
			\item[(a)] $\tilde{{\mathcal F}}(\la, w) \neq w$ for all $\la \in [0, \hat{\la}]$ whenever $\|w\|_{C(\overline\Om)} \in \{r, R\}$; and
			\item[(b)] $
			\text{deg}(I-\tilde{{\mathcal F}}(\la, \cdot), B_R \setminus \overline{B}_r, 0)=-1
			$ for all $\la \in [0, \hat{\la}]$.
		\end{description}
		
	\end{lem}
	
	\begin{proof}
		
		(a) Suppose not. Then there exist sequences $\la_n \geq 0$ with $\la_n \to 0$ and $w_n \in C(\overline \Om)$ such that 
		$\tilde{{\mathcal F}}(\la_n, w_n) = w_n$ and $\|w_n\|_{C(\overline\Om)}=r$ (or $\|w_n\|_{C(\overline\Om)}=R$). Since $w_n$ is bounded and $\tilde{\mathcal{F}}$ is compact, $(\la_n, w_n) \to (0, w)$ for some $w \in C(\overline \Om)$ with $\|w\|_{C(\overline\Om)}=r$ or $\|w\|_{C(\overline\Om)}=R$, a contradiction to Lemma~\ref{lem:0:r} or Lemma~\ref{lem:0:R}, respectively. Hence there exists $\hat{\la}>0$ satisfying (a).
		\par (b) Now using $\la \in [0, \hat{\la}]$ as the homotopy parameter, it follows from part (a) that 
		$$
		\text{deg}(I-\tilde{{\mathcal F}}(\la, \cdot), B_R \setminus \overline{B}_r, 0)=\text{const.} \quad \forall \la \in [0, \hat{\la}]\,.
		$$
		In particular, it follows from \eqref{deg:0:-1}  that  for all $\la \in [0, \hat{\la}]$
		$$
		\text{deg}(I-\tilde{{\mathcal F}}(\la, \cdot), B_R \setminus \overline{B}_r, 0)= 
		\text{deg}(I-\tilde{{\mathcal F}}(0, \cdot), B_R \setminus \overline{B}_r, 0)=-1\,.
		$$
		This complete the proof of Lemma~\ref{lem:rescaled}.
	\end{proof}
	
	Lemma~\ref{lem:rescaled} implies that the re-scaled problem \eqref{pde:rescaled} has a nontrivial solution $w_{\la} \in C(\tilde \Om)$ for all $\la \in [0, \hat{\la}]$ satisfying $r < \|w_{\la}\|_{C(\tilde\Om)} < R$. Moreover, since $f$ is nonnegative and satisfies  $\mathrm{(H)}_\infty$, so does $\tilde{f}$ and hence 
	$w_{\la}>0$ in $\overline{\Om}$ by Proposition~\ref{pro:v>0}.
	\par Now we return to the original problem \eqref{pde}. Using the re-scaling  
	$$
	u=\la^{-\frac{1}{p-1}}w_{\la}\,,
	$$
	we can conclude that \eqref{pde} has a positive solution $(\la, u)$ for $\la \in (0, \hat{\la}]$. Also, since $\|w_{\la}\|_{C(\overline\Om)}>r>0$, it follows that 
	$\|u\|_{C(\overline\Om)} \to +\infty$ as $\la \to 0^+.$
	\par We use the following Leray-Schauder continuation theorem  to establish the last part of Theorem~\ref{thm:main}.
	\begin{pro}(\cite[Prop.~2.3]{deF-Lion_Nus_1982})
		\label{prop:continuation}
		Let $X$ be a Banach space and $U$ a bounded open subset of $X$. Let $\mathcal{T}: [a, b] \times \overline{U} \to X$ be a compact map and $\mathcal{S}:= \left\{(t, x) \in [a, b] \times U:\, \mathcal{T}(t, x)=x\right\}$ is the set of all fixed points of $\mathcal{T}$. Assume that  
		\begin{itemize}
			\item $\mathcal{T}(t, x) \neq x$ for $(t, x) \in [a, b] \times \p U$;
			\item $\text{deg}(I - \mathcal{T}(t, \cdot), 0) \neq 0$ for all $t \in [a, b]$.
		\end{itemize}
		Then there exists a connected component $\mathcal{D}$ of $\mathcal{S}$ such that $\mathcal{D} \cap (\{a\}\times U)$ and  $\mathcal{D} \cap (\{b\}\times U)$ are nonempty.
	\end{pro}
	Now, by taking  $[a, b]=[0, \hat{\la}]$, $U=B_R \setminus \overline{B_r}$, $\mathcal{T}=\tilde{{\mathcal F}}(\cdot, \cdot)$ in Proposition~\ref{prop:continuation},
	it follows, using 
	Lemma~\ref{lem:rescaled}, that the re-scaled problem \eqref{pde:rescaled} has a connected component $\mathcal{D}$ of positive weak solutions along which $\la$ takes all values in $[0, \hat{\la}]$. This in turn, again using $
	u=\la^{-\frac{1}{p-1}}w_{\la}
	$, implies that there exists a connected component $\mathscr{C}^+ \subset \Sigma$ of positive weak solutions of  \eqref{pde} bifurcating from infinity at $\la_{\infty}=0$. This completes the proof of Theorem~\ref{thm:main}.
	\section{Global Bifurcation}
	\label{sec:proof:thm:bif}
	\smallskip
	In this section, we will prove that there exists a connected set of positive weak  solutions $\mathscr{C}^{+}$ of \eqref{pde}  bifurcating from the trivial solution  at $\la=\frac{\mu_1}{f'(0)}$, and bifurcating from infinity at $\la =0$.
	Furthermore, we discuss the direction of bifurcation 
	of positive weak solutions at $(\frac{\mu_1}{f'(0)}, 0)$.
	Finally, we prove Theorem~\ref{th:conexion:0:inf}.
	
	\medskip
	\subsection{Bifurcation from the trivial solution}
	We  first show that the condition $\mathrm{(H)}_0$ guarantees solutions bifurcating from the trivial solution. The proof is similar to the case of bifurcation from infinity, see for instance \cite[Proposition 3.1]{Arrieta-Pardo-RBernal_2007}. We provide the proof below for completeness. 
	
	\begin{pro}\label{pro:convergence}
		Assume that  the nonlinearity  $f\in C^1([0,\infty))$   
		satisfies the hypothesis $\mathrm{(H)}_0$. 
		Let $\{\la_n\}$ be a convergent sequence of real numbers and $u_n$  be the corresponding sequence of positive weak solutions of equation  \eqref{pde} satisfying  $||u_n||_{C(\overline{\Om})}\to 0$ as $n\to  \infty$. Then, necessarily  $\la_n\to\dfrac{\mu_1}{f'(0)}$,
		and   $\{u_{n}\}$ satisfies, up to a subsequence,  
		\begin{equation*}\label{eq:conv:Phi1}
			\displaystyle \frac{u_{n}}{||u_{n}||_{C(\overline{\Om})}}\to \varphi_1\;\;\;{\text{in}}\;C^{\be  }(\overline{\Om}) 
		\end{equation*}
		for some $\be  \in (0, 1)$.
	\end{pro}
	\begin{proof}  Suppose that $\la_n \to \underline{\la}$ for some $\underline{\la} \in \R$ and set $\displaystyle v_n:=\frac{u_{n}}{||u_{n}||_{C(\overline{\Om})}}$.
		Observe that $v_{n}$ is a weak solution of the problem 
		\begin{equation}\label{eqU1}
			\begin{cases}
				\displaystyle
				-\Delta v_{n} + v_{n} = 0 &\text{in}\; \Om\,;\\
				\displaystyle\frac{\p  v_{n}}{\p  \eta}=  \la_{n}  f'(0) v_n
				+\la_{n}\frac{\mathcal{R}(u_{n}) }{||u_{n}||_{C(\overline{\Om})}}
				& \text{on}\;
				\p \Om\,.
			\end{cases}
		\end{equation}
		It follows from $\mathrm{(H)}_0$  that  
		$
		\frac{\mathcal{R}(u_{n})}{||u_{n}||_{C(\overline{\Om})}}\to 0$ in $C(\overline{\Om}) $ as $n\to \infty$. Therefore, the right-hand side of the second equation in \eqref{eqU1} is bounded in $C(\overline \Om)$. Hence, by the elliptic regularity, $v_n\in W^{1,s}(\Om)$ for any  $s>1$, in particular for $s>N$. Then, the Sobolev embedding theorem implies that $||v_n||_{C^{\al }(\overline{\Om})}$ is bounded by a constant $C$ that is independent of $n$.
		Then, the compact embedding of $C^{\al }(\overline{\Om})$ into $C^{\be  }(\overline{\Om})$ for $0<\be  <\al $  yields, up to a subsequence,  $v_{n}\to \Phi \geq 0$ in  $C^{\be  }(\overline{\Om})$. Since $||v_{n}||_{C(\overline{\Om})}=1$, we have that $||\Phi||_{C(\overline{\Om})}=1$. Hence, $\Phi\not\equiv 0.$
		
		\par Using the weak formulation of equation (\ref{eqU1}),  passing to the limit, and taking into account  that $\la_n \to \underline{\la}$ for some  $\underline{\la}\in\mathbb{R}$ and $v_{n}\to \Phi$, we obtain that  $\Phi$ is a weak solution of the equation
		\begin{equation*}
			\begin{cases}
				\displaystyle
				-\Delta \Phi + \Phi = 0 &\text{in}\; \Om\,;\\
				\displaystyle\frac{\p  \Phi}{\p  \eta}=  \underline{\la}  f'(0)\Phi & \text{on}\;
				\p \Om.
			\end{cases}
		\end{equation*}
		Then, it follows that  $\underline{\la}f'(0)=\mu_1$,  the first Steklov eigenvalue, and $\Phi=\varphi_1>0$ is its corresponding eigenfunction, ending the proof. 
	\end{proof}
	
	Now, we will show that  $\big(\frac{\mu_1}{f'(0)},0\big)$ is a bifurcation point from the trivial solution of positive weak solutions of \eqref{pde}. That is, there exists a sequence $(\la_n,u_n) \in \Sigma$   such that $\la_n\to \frac{\mu_1}{f'(0)}$, $u_n>0$ on $\overline{\Om}$, and that $||u_n||_{C(\overline{\Om})}\to 0.$
	In particular, we have the following result.
	
	\begin{thm}\label{th:bif:0}
		Assume that  the nonlinearity  $f\in C^1([0,\infty))$  
		satisfies hypothesis $\mathrm{(H)}_0$.   Then, there exists a connected component ${\mathscr
			C}^+ \subset \Sigma$ of positive weak solutions of \eqref{pde} emanating from the trivial solution at 
		$\big(\frac{\mu_1}{f'(0)},0\big)\in\R\times C(\overline{\Om})$.
		Moreover, ${\mathscr
			C}^+$ 
		is  unbounded in $\R\times C(\overline{\Om})$. 
	\end{thm}
	
	\begin{proof} The proof follows from the general results on bifurcation from the trivial solutions given in \cite[Thm.~1.3]{R71}. 
		More precisely, there exists a connected component ${\mathscr
			C}^+ \subset \Sigma$ of positive weak solutions of \eqref{pde} emanating from the trivial solution at 
		$\big(\frac{\mu_1}{f'(0)},0\big)\in\R\times C(\overline{\Om})$ and,  the branch ${\mathscr C}^+ $ either meets another bifurcation point from the trivial solution, or it is unbounded in $\R\times C(\overline{\Om})$.
		Since $f\ge 0$ satisfies  $\mathrm{(H)}_0$, it follows from  Lemma \ref{lem:reg:lbp}~(iv) and Proposition~\ref{pro:v>0} that the branch contains only positive solutions. 
		From the Crandall-Rabinowitz Theorem, see \cite{C-R}, ${\mathscr C}^+ $  can neither meet another bifurcation point from zero (that is, another point
		$\big(\frac{\mu'}{f'(0)},0\big)$ for another Steklov eigenvalue $\mu'$), nor can meet $\big(\frac{\mu_1}{f'(0)},0\big)$
		again, so the branch is unbounded in $\R\times C(\overline{\Om})$. 
	\end{proof} 
	\subsection{Subcritical and supercritical bifurcations from the trivial solution}
	
	In this subsection, we discuss sufficient conditions for the bifurcation from the trivial solution to be either subcritical (to the left) or supercritical (to the right). Following lemma is key in determining the direction of bifurcation from the trivial solution at $\left(\frac{\mu_1}{f'(0)}, 0\right)$.
	
	\begin{lem}
		Assume that  the nonlinearity  $f\in C^1([0,\infty))$   
		satisfies the hypothesis $\mathrm{(H)}_0$.  Consider a sequence of positive weak solutions $u_n$ of \eqref{pde} corresponding to the parameters $\la_n$ such that $\la_n\to \frac{\mu_1}{f'(0)}$ and
		$\|u_n\|_{C(\overline{\Om})}\to 0
		$. 
		Then,  we have 
		\begin{align}\label{ineq-positive-sol}
			\underline{\mathcal R}_0\ \frac{\mu_1}{\big(f'(0)\big)^2}\frac{\int_{\p\Om}\varphi_1^{1+\nu}}{
				\int_{\p\Om} \varphi_1^2}
			&\leq \liminf_{n\to\infty}\frac{\frac{\mu_1} {f'(0)}-\la_n}{
				\|u_n\|_{C(\overline{\Om})}^{\nu -1}}\\
			&\leq \limsup_{n\to\infty}\frac{\frac{\mu_1} {f'(0)}-\la_n}{
				\|u_n\|_{C(\overline{\Om})}^{\nu-1}}
			\leq \overline{\mathcal R}_0\ \frac{\mu_1}{\big(f'(0)\big)^2}\frac{\int_{\p\Om} \varphi_1^{1+\nu}}{\int_{\p\Om} \varphi_1^2}\,,\nonumber
		\end{align}
		where $\underline{\mathcal R}_0$ and $\overline{\mathcal R}_0$ are defined in \eqref{R}, and $\nu>1$ as defined in $\mathrm{(H)}_0$.
	\end{lem}
	
	\begin{proof}
		Using the weak formulation of \eqref{pde} with $\varphi_1$ as the test function, we get
		$$
		\int_{\Om}\nabla u_n \nabla \varphi_1 + \int_{\Om}u_n\varphi_1 = \la_n f'(0)\int_{\p \Om}u_n\varphi_1+\la_n
		\int_{\p \Om}\mathcal{R}(u_n)\varphi_1\,,
		$$
		which yields
		$$
		(\mu_1-\la_n f'(0))\int_{\p \Om}u_n\varphi_1 = \la_n
		\int_{\p \Om}\mathcal{R}(u_n)\varphi_1\,.
		$$
		Consequently, we get
		\begin{equation}
			\label{eq:bdry0:1}
			\frac{(\mu_1-\la_n f'(0))} {||u_n||_{C(\overline{\Om})}^{\nu-1}} \int_{\p \Om}\frac{u_n}{||u_n||_{C(\overline{\Om})}}\varphi_1 = 
			\la_n \int_{\p \Om}\frac{\mathcal{R}(u_n)}{||u_n||_{C(\overline{\Om})}^{\nu}}\varphi_1.
		\end{equation}
		From Fatou's Lemma,
		\begin{align}
			\label{eq:bdry0:2}
			\displaystyle\liminf_{n\to \infty} & \displaystyle\int_{\p \Om }\frac{\mathcal{R}(u_n)}{ u_n^\nu }\left(\frac{u_n}{ \|u_n\|_{C(\overline{\Om})}}\right)^\nu\varphi_1 \\
			&\geq  \displaystyle\int_{\p \Om } \liminf_{n\to \infty}\left[\frac{\mathcal{R}(u_n)}{ u_n^\nu }\left(\frac{u_n}{ \|u_n\|_{C(\overline{\Om})}}\right)^\nu\varphi_1\right]  \geq  \displaystyle \underline{\mathcal R}_0 \int_{\p \Om } 
			\varphi_1^{1+\nu }\,,\nonumber
		\end{align}
		where we have used the definition of $\underline{\mathcal R}_0$ (see \eqref{R}), that $ \varphi_1>0$ on $\p\Om $ and the fact that $\frac{u_n}{\|u_n\|_{C(\overline{\Om})}}\to \varphi_1$ uniformly on $\p\Om$ ( see Proposition \ref{pro:convergence}). 
		
		Passing to the limit in \eqref{eq:bdry0:1} and using \eqref{eq:bdry0:2}, we obtain the first inequality of \eqref{ineq-positive-sol}. The second inequality is trivial and the third is obtained likewise. 
	\end{proof}

	Now, we can state the following result with regarding subcritical or supercritical bifurcations from the trivial solution. 
	\begin{thm} ({\bf Bifurcation of positive solutions from the trivial solution}) 
		\label{th:bif-from-first-sub}
		Assume that  the nonlinearity  $f\in C^1([0,\infty))$   
		satisfies hypothesis $\mathrm{(H)}_0$.  
		Then, the following holds.
		
		\begin{itemize}
			\item[\rm (i)]{\bf (Subcritical bifurcations)}. If
			$
			\underline{\mathcal R}_0
			> 0 
			$,
			then the bifurcation of positive  weak solutions from the trivial solution 
			at $\la  =\frac{\mu_1}{f'(0)}$ is subcritical, i.e. $\la  < \frac{\mu_1}{f'(0)} $ for every positive  solution $(\la  , u)$ of
			\eqref{pde} with $(\la ,\|u\|_{C(\overline{\Omega})})$ in a neighborhood of
			$(\frac{\mu_1}{f'(0)}, 0)$.
			\item[\rm (ii)]{\bf (Supercritical bifurcations)}. If
			$
			\overline{\mathcal R}_0< 0 
			$,
			then the bifurcation of positive  weak solutions from the trivial solution
			at $\la  =\frac{\mu_1}{f'(0)}$ is supercritical, i.e. $\la  > \frac{\mu_1}{f'(0)} $
			for every positive  solution $(\la ,u)$ of
			\eqref{pde} with $(\la , \|u\|_{C(\overline{\Omega})})$ in a neighborhood of
			$(\frac{\mu_1}{f'(0)},0)$.
		\end{itemize}
	\end{thm}
	
	\begin{proof}
		Consider a sequence of positive weak solutions $u_n$ of \eqref{pde} corresponding to the parameters $\la_n$ such that $\la_n\to \frac{\mu_1}{f'(0)}$ and
		$\|u_n\|_{C(\overline{\Om})}\to 0
		$.  Observe that, by \eqref{ineq-positive-sol}, conditions
		$\underline{\mathcal R}_0>0$ and $\overline{\mathcal R}_0<0$ imply that $\frac{\mu_1}{f'(0)}>\la _n$ and $\frac{\mu_1}{f'(0)}<\la _n$, respectively, for sufficiently large $n$. This completes the proof.
	\end{proof}
	\subsection{Proof of Theorem \ref{th:conexion:0:inf}}
	The proof will be completed in several steps.
	
	\smallskip
	\noindent{\bf Step~1:} By Theorem \ref{th:bif:0}, there exists a connected component $\mathscr{C}^+$ of positive weak solutions of \eqref{pde}  bifurcating from the trivial solution at the bifurcation point
	$\big(\frac{\mu_1}{f'(0)},0\big)$ and that $\mathscr{C}^+$  is unbounded in $\R\times C(\overline{\Om})$.
	
	\smallskip
	\noindent{\bf Step~2:} 
	At this step, we show that  \eqref{pde} has no positive weak solution for $\la  > \frac{\mu_1}{K}$, where $K>0$ is as given in the hypothesis \eqref{hyp:nonexist}.

	Indeed, let $u $ be a positive weak solution of \eqref{pde} for some $\la >0$. Then, using $\varphi_1\ge 0$ as the test function, we get 
	\begin{align*}
		0 &= \int_{\Om}\nabla u \nabla \varphi_1 + \int_{\Om}u\varphi_1 -  \la  \int_{\p \Om}f(u)\varphi_1\\
		&= \la  \int_{\p \Om}f(u)\varphi_1 - \int_{\Om}\nabla u \nabla \varphi_1 - \int_{\Om}u\varphi_1\\
		&\geq   \la  K\int_{\p \Om} u \varphi_1 - \int_{\Om}\nabla u \nabla \varphi_1 - \int_{\Om}u\varphi_1  
		= (\la  K- \mu_1)\int_{\p \Om} u\varphi_1\,.
	\end{align*}
	This yields $\la  \leq \frac{\mu_1}{K}$. Hence there exists no positive weak solution $u$ of \eqref{pde} for $\la  > \frac{\mu_1}{K}$,
	completing the proof of this step.
	
	\smallskip
	
	\noindent{\bf Step~3:} Here, we show  that $\mathscr{C}^+$ from Step~1 contains weak positive solutions that bifurcate from infinity at $\lambda=0$, and establish  \eqref{unbdd:zero}. 
	\par By Step~1-Step~2, if $(\la, u) \in \mathscr{C}^+$ then 
	$
	\|u\|_{C(\overline{\Om})} \to 0 \mbox{ as } \la\to\frac{\mu_1}{f'(0)},
	$ 
	and $\mathscr{C}^+$ is bounded in the $\la$-direction. Hence, 
	there exists a sequence  $(\la_n, u_n) \in \mathscr{C}^+$ such that  $\la_n \in (0, K)$ and $\|u_n\|_{C(\overline{\Om})}\to \infty$.
	By choosing a subsequence if necessary, there exists a sequence  $(\la_n, u_n) \in \mathscr{C}^+$ with the property that $\la_n\to \tilde{\la}$ and $\|u_n\|_{C(\overline{\Om})}\to \infty$. It suffices to show $\tilde{\la}=0$. 
	\par Assume to the contrary that $\tilde{\la}> 0$.  For $a_0>0,$ let  $[a_0,b_0]$ be  any  fixed compact interval with $\tilde{\la}\in(a_0,b_0)$. By   Proposition \ref{prop:apriori}, for any $\la\in[a_0,b_0]$,   there exists a uniform constant $M=M(a_0,b_0) > 0$ such that  for every $(\la,w)$ with $\la\in[a_0,b_0]$ and $w$  a positive  weak solution  of the re-scaled problem \eqref{pde:extend:2}$_\la$,  we have
	\begin{equation*}
		\|w\|_{C(\overline{\Om})}\leq M\,.
	\end{equation*}
	
	Here we recall from Section~\ref{sec:proof:thm:main} that for  any $\la>0$, $u$ is a positive weak solution of \eqref{pde} if and only if $w=\la ^{\frac{1}{p-1}}u$ is a weak solution of  \eqref{pde:extend:2}. Hence, 
	\begin{equation}
		\label{norm:w:to:u}
		\|u\|_{C(\overline{\Om})}
		\leq \la^{-\frac{1}{p-1}}M\leq a_0^{-\frac{1}{p-1}}M=:M'\qq{for  any} \la\in[a_0,b_0],
	\end{equation} 
	which contradicts that $\|u_n\|_{C(\overline{\Om})}\to \infty$ with $\la_n\to \tilde{\la}>0$. Hence $\tilde{\la}=0$. As a conclusion, necessarily, $\mathscr{C}^+$ contains a unique bifurcation point from infinity at $\la=0$ and \eqref{unbdd:zero} holds. Then,  \eqref{pde} has a positive weak solution for any $\lambda \in \big(0, \frac{\mu_1}{f'(0)}\big)$. This completes Step~3.
	\smallskip
	Now, set 
	\begin{equation*}
		\label{lambda:sup}
		\bar{\la}:=\sup\{\la>0: (\la, u) \in \mathscr{C}^+ \}.
	\end{equation*}
	Then, $\bar{\la}<\infty$ by Step~2.
	\smallskip
	\noindent{\bf Step~4:} Assuming   $\overline{{\mathcal R}}_0<0$, we prove the existence of two positive weak solutions for each $\la \in \left(\frac{\mu_1}{f'(0)}, \bar\la\right)$. 
	\par It follows from Theorem \ref{th:bif-from-first-sub} (ii), that the bifurcation is  supercritical at the bifurcation point  $\big(\frac{\mu_1}{f'(0)},0\big)$ from the  trivial solution. 
	Note that since   $\overline{{\mathcal R}}_0<0$, $\bar{\la}> \frac{\mu_1}{f'(0)}$. 
	Let $\la_0 \in \left(\frac{\mu_1}{f'(0)}, \overline{\la}\right)$ and $u_0$ be a positive weak solution corresponding to $\la_0$. Now, let $\la \in \left(\frac{\mu_1}{f'(0)}, \la_0\right)$ be fixed. We show that there exist two distinct positive weak solutions of \eqref{pde} corresponding to $\la$ using degree theory. For this, first we extend $f$ to $\mathbb{R}$ by setting $f(t)=0$ for $t < 0$. \\
	\smallskip
	\noindent{\bf First solution corresponding to $\la$:} 
	\par  First we note that, since $f$ is Lipschitz continuous, there exists $c \in \mathbb{R}$ such that $\lambda f(s) + cs$ is nondecreasing on $[0, M']$, 
	where $M' > M $, and $M>0$ is  given by Proposition~\ref{prop:apriori}. Now let  $\theta \in [0, 1]$ and $\beta > \mu_1$. For a given $u \in C(\overline \Omega)$, define the operator $T_{\theta}: C(\overline \Omega) \to C(\overline \Omega)$ by  $v=T_{\theta}(u):=(S\circ f_{\theta}\circ\Gamma)u$, where $v$ is given by
	\begin{equation*}
		\label{pde:theta:homotopy}
		\left\{
		\begin{array}{rcll}
			-\Delta v +v &=&  0 \quad \mbox{in}\quad \Om\,;\\
			\frac{\p v}{\p \eta} + \theta c\,v&=& \theta (\lambda f(u) + cu) + (1-\theta)(\beta u^++1)\quad \mbox{on}\quad \p \Om\,,
		\end{array}
		\right.
	\end{equation*}
	and $f_{\theta}(u):=\theta (\lambda f(u) + cu) + (1-\theta)(\beta u^++1)$. We note that 
	$T_{\theta}$
	is compact by Remark~\ref{rem:trace:u:reg}, and fixed point of the operator $T_1$ is a weak solution of \eqref{pde}.
	\par We begin by establishing that $u_0 > \epsilon \varphi_1$ for sufficiently small $\epsilon >0$.
	Clearly,  $u_0-\epsilon\varphi_1$ satisfies
	\begin{equation*}
		\label{u:ephi_1}
		-\Delta(u_0-\epsilon \varphi_1) + (u_0-\epsilon \varphi_1) = 0 \mbox{ in }  \Om\,.
	\end{equation*}
	Now, using the hypothesis \eqref{hyp:nonexist}, and the facts that $\la > \frac{\mu_1}{f'(0)}$, $\|u_0\|_{C(\overline \Omega)} < M'$ and $f$ is continuous,  we get
	$$
	\la f(u_0) - \epsilon \mu_1 \varphi_1 \geq \frac{\mu_1}{f'(0)} \left( f(u_0) - \epsilon f'(0)\varphi_1\right) \geq 0
	$$
	for  $\epsilon>0$ sufficiently small. Then
	\begin{equation*}
		\label{ineq:compare:u:ephi_1}
		\frac{\partial (u_0-\epsilon \varphi_1) }{\partial \eta} =
		\lambda f(u_0)  - \epsilon \mu_1 \varphi_1 \geq 0 \mbox{ on } \partial \Om\,.
	\end{equation*}
	Therefore, by Proposition~\ref{pro:v>0}, $u_0 > \epsilon \varphi_1 $ for $\epsilon>0$ sufficiently small. \\

	Now define
	$$
	Y:= \left\{v \in C(\overline \Om):    \|v\|_{C(\overline \Om)} < M' \mbox{ and  }  v > \epsilon \varphi_1 \mbox{ on } \overline \Om
	\right\}\,,
	$$
	and 
	$$
	Z:=\left\{v \in Y: v < u_0 \mbox{ on }  \overline{\Om} \right\}\,,
	$$
	where  $\epsilon>0$ to be chosen sufficiently small later such that in particular $u_0>\epsilon \varphi_1$ in $\Omb$.\\ 
	
	\noindent{\bf Claim I:} $\text{deg}(I-T_1, Y, 0)=0$.
	\par First, we justify that the degree $\text{deg}(I-T_{\theta}, Y, 0)$ is well defined and independent of  $\theta \in[0, 1]$. That is, $u \neq T_{\theta}u$ for any $u$ on the boundary of $Y$, $\partial Y$. We note that if $u \in \partial Y$, then either $\|u\|_{C(\overline \Omega)}=M'$ or $u=\epsilon \varphi_1$. Now, if $\|u\|_{C(\overline \Omega)}=M'$, then
	by Proposition~\ref{prop:apriori},  $u \neq T_{\theta}u$ for any $\theta \in [0, 1]$. On the other hand, if  $u=\epsilon \varphi_1$ is a solution  of $u=T_{\theta}u$ for  $\theta =0$, then $\beta > \mu_1$ yields 
	the contradiction
	$$
	\beta \epsilon\varphi_1 +1 =\frac{\partial \epsilon \varphi_1}{\partial \eta}=\mu_1 \epsilon \varphi_1< \beta \epsilon\varphi_1\,.
	$$
	Thus, $u \neq T_{\theta}u$ when $u=\epsilon \varphi_1$.
	
	\par Now, by repeating arguments in Step~2 with $\lambda f(u)$ replaced by $\beta u^+ + 1$ and using $\beta > \mu_1$, we see that $u \neq T_{0}u$ for any $u \in Y$. Then, using $\theta \in [0, 1]$ as a homotopy parameter, we conclude that 
	\begin{equation}
		\label{deg:0}
		\text{deg}(I-T_1, Y, 0)=\text{deg}(I-T_{\theta}, Y, 0)=\text{deg}(I-T_{0}, Y, 0)=0\,.
	\end{equation}
	\noindent{\bf Claim II:} $\text{deg}(I-T_1, Z, 0)=1$.\\
	We fix $\psi_0 \in Z$ and show $\text{deg}(I-(\theta T_1 + (1-\theta)T_{\psi_0}), Z, 0)=1$ for $\theta \in [0, 1]$, where   $T_{\psi_0}$ maps every element of $Z$ to $\psi_0$. By 
	$v=(\theta T_1 + (1-\theta)T_{\psi_0})u$, for $\theta \in [0, 1]$, we mean 
	\begin{equation*}
		\label{pde:theta:homotopy:2}
		\left\{
		\begin{array}{rcll}
			-\Delta v +v &=&  0 \quad \mbox{in}\quad \Om\,;\\
			\frac{\p v}{\p \eta} + \theta c\,v&=& \theta (\lambda f(u) + cu) + (1-\theta)\psi_0 \quad \mbox{on}\quad \p \Om\,.
		\end{array}
		\right.
	\end{equation*}
	Now we show that $\text{deg}(I-(\theta T_1 + (1-\theta)T_{\psi_0}), Z, 0)$ is well defined and  independent of $\theta \in [0, 1]$. 
	Indeed,  note that if $u \in \overline{Z}$, that is, $u \leq u_0$, then by Proposition~\ref{pro:v>0} $v=T_1u \in Z$, since $-\Delta v + v=0$ in $\Omega$ and
	$$
	\frac{\partial v}{\partial \eta} + cv = \lambda f(u) + cu \leq  \la f(u_0) + cu_0 < \la_0 f(u_0) + cu_0 \quad \mbox{ on } \partial \Omega\,.
	$$
	Also, $T_{\psi_0}u \in Z$ for $u \in Z$. Then $\theta T_1u + (1-\theta)T_{\psi_0}u  \in Z$ for all $\theta \in [0, 1]$, since $Z$ is convex. Hence, there is no solution of $I- (\theta T_1 + (1-\theta)T_{\psi_0})$ on the boundary of $Z$, and 
	$\text{deg}(I- (\theta T_1 + (1-\theta)T_{\psi_0}), Z, 0)$ is well defined for all $\theta \in [0, 1]$. Therefore, since $\psi_0 \in Z$, we have
	\begin{equation}
		\label{deg:1}
		\text{deg}(I- T_1, Z, 0) = \text{deg}(I- T_{\psi_0}, Z, 0) = \text{deg}(I, Z, \psi_0)=1\,,
	\end{equation}
	completing Claim II.
	\par Combining \eqref{deg:0} and \eqref{deg:1}, one has
	$\text{deg}(I- T_1, Y \setminus \overline{Z}, 0)=-1$ and hence there exists a positive weak solution $u_2 \in Y \setminus \overline{Z}$ of \eqref{pde} corresponding to the fixed $\la$.\\ \smallskip
	\noindent{\bf Second solution corresponding to $\la$:}  
	\par We construct the second positive weak solution distinct from $u_2$ by the method of sub- and supersolutions. Using the facts that $f(0)=0$ and $f'(0)>0$, we verify that $\underline{u}=\epsilon \varphi_1$ is a subsolution of \eqref{pde} for $\epsilon \approx 0$. Indeed, we observe that since $\lambda > \frac{\mu_1}{f'(0)}$ is fixed, $\xi(s):=\mu_1 s -\lambda f(s)$ satisfies $\xi(0) = 0$ and $\xi'(0) < 0$, then $\xi(s) < 0$ for $s \approx 0$. Therefore, for all $0 \leq \psi \in H^1(\Om)$, the following holds for $\epsilon \approx 0$
	$$
	\int_{\Om}\nabla \underline{u} \nabla \psi + \int_{\Om}\underline{u}\psi =\mu_1 \int_{\p \Om}(\epsilon \varphi_1)\psi\leq   \la \int_{\p \Om}f(\epsilon \varphi_1)\psi=\la \int_{\p \Om}f(\underline{u})\psi\,. 
	$$
	Note that $u_0 \in Y$ since $\epsilon \varphi_1 < u_0 < M< M' $ for sufficiently small $\epsilon>0$.
	It follows from \cite{BCDMP-2} that $\text{min}(u_2, u_0)$ is a strict supersolution of \eqref{pde}.
	Since $u_0, u_2 \in Y$,  $\underline{u}=\epsilon \varphi_1 < \text{min}(u_2, u_0)$ on $\overline{\Omega}$. Hence, there exists a positive weak solution  $u_1$ of \eqref{pde} corresponding to the fixed $\la$ satisfying 
	$\epsilon \varphi_1 \leq u_1 < u_2$ on $\Omega$  by Proposition~\ref{pro:v>0}. This completes Step~4.
	
	\smallskip
	\noindent{\bf Step~5:} At this step,  we prove the existence of a solution for $\la=\overline{\la}$. For each $\la\in (\frac{\mu_1}{f'(0)},\bar{\la})$, 
	problem \eqref{pde} admits a 
	positive weak solution $u_\la$. 
	\par Using Proposition~\ref{prop:apriori}, \eqref{norm:w:to:u} for $\lambda \in [\frac{\mu_1}{f'(0)},\bar{\la}]$, and Proposition~\ref{pro:0:inf},  there exists a uniform constant $C>0$ such that 
	$\|u_{\la}\|_{C^\al(\Omb)} \le C$ 
	for any $\la\in (\frac{\mu_1}{f'(0)},\bar{\la})$. By compact embeddings, $u_{\la}$ has a  subsequence that converges to (say), $u_{\, \bar{\la}}$  in $C^{\be}(\Omb )$
	as $\la \to \bar\la$,  where $\be<\al$. 
	\par Moreover,
	\begin{equation*}\label{very:weak:sol:2}
		\|u_{\la}\|_{H^1(\Om)}^2 =
		\int_{\Om}  |\na u_{\, \la} | ^2 + \int_{\Om}|u_{\, \la}|^2 = \la\int_{\p \Om}f(u_{\, \la})u_{\, \la}\le C ,\quad \forall \la \in \left(\frac{\mu_1}{f'(0)},\bar{\la}\right)\,.
	\end{equation*}
	By the reflexivity of ${H^1(\Om)}$, $u_{\la}$ has a  subsequence that converges weakly to (say), $u_{\, \bar{\la}}$  in  ${H^1(\Om)}$ as $\la\to \bar\la$. 
	On the other hand, since $u_{\, {\la}}\to u_{\, \bar{\la}}\in C^{\be}(\Omb )$ and $f$ is locally Lipschitz, then $f(u_{\, {\la}})\to f(u_{\, \bar{\la}})$  in $C^{\be}(\Omb )$ as  $\la \to\bar \la$. 
	\par Then, by taking limits in the weak formulation of $u_{\la}$ as $\la \to \bar\la$, we get 
	\begin{equation*}\label{very:weak:sol:3}
		\int_{\Om}  \na u_{\, \bar{\la}} \na \psi + \int_{\Om}u_{\, \bar{\la}}\psi = \bar{\la}\int_{\p \Om}f(u_{\, \bar{\la}})\psi\,.
	\end{equation*}
	Hence $u_{\, \bar{\la}}$ is a positive  weak solution of \eqref{pde}$_{\bar{\la}}$.
	
	\par Therefore, \eqref{pde} has at least two  positive weak solutions for $\la \in \big(\frac{\mu_1}{f'(0)}, {\bar \la}\big)$, and at least one  positive weak solution for $\la =\bar \la$. Finally, since the connected set $\mathscr{C^+}$ bifurcates to the right at $\big(\frac{\mu_1}{f'(0)}, 0\big)$ and bifurcates from infinity at $\la=0$, $\mathscr{C^+}$ must cross the hyperplane $\la =\frac{\mu_1}{f'(0)}$ at a point distinct from $u=0$. Hence, the problem \eqref{pde} has a positive weak solution for
	$\la =\frac{\mu_1}{f'(0)}$. This completes  the proof of Theorem~\ref{th:conexion:0:inf}. \hfill $\Box$
	
	\medskip
	\noindent{\bf Acknowledgements:}
	This project was approved by MSRI  Summer Research for  Women in
	Mathematics (SWiM) Program for summer 2020. The visit was postponed due to Covid-19, but authors are grateful to MSRI for bringing together this group for collaboration.
	R. Pardo was partially supported by   grant PID2019-103860GB-I00,  MICINN,  Spain, and by UCM-BSCH, Spain, GR58/08, Grupo 920894.

\end{document}